\documentclass[journal,twoside,web]{ieeeconf}

\usepackage{cite}
\usepackage{amsmath,amssymb,amsfonts}
\usepackage{algorithmic}
\usepackage{graphicx}
\usepackage{textcomp}
\usepackage{comment}
\usepackage{subcaption}

\usepackage[font=small,labelfont=bf]{caption}

\usepackage[colorlinks = true,
            linkcolor = blue,
            urlcolor  = blue,
            citecolor = blue,
            anchorcolor = blue]{hyperref}

\newtheorem{theorem}{Theorem}

\newtheorem{proposition}{Proposition}

\newtheorem{lemma}{Lemma}[section]

\newtheorem{assumption}{Assumption}
\newtheorem{problem}{\textbf{Problem}}

\newtheorem{remark}{Remark}

\usepackage[nameinlink]{cleveref}
\crefname{equation}{}{}
\crefname{theorem}{Theorem}{Theorems}
\crefname{corollary}{Corollary}{Corollaries}
\crefname{example}{Example}{Examples}
\crefname{assumption}{Assumption}{Assumptions}
\crefname{lemma}{Lemma}{Lemmas}
\crefname{proposition}{Proposition}{Propositions}
\crefname{figure}{Figure}{Figures}
\crefname{table}{Table}{Tables}
\crefname{fact}{Fact}{Facts}
\crefname{conjecture}{Conjecture}{Conjectures}
\crefname{section}{Section}{Sections}
\crefname{appendix}{Appendix}{Appendices}
\Crefname{equation}{}{}
\Crefname{theorem}{Theorem}{Theorems}
\Crefname{corollary}{Corollary}{Corollaries}
\Crefname{example}{Example}{Examples}
\Crefname{lemma}{Lemma}{Lemma}
\Crefname{proposition}{Proposition}{Proposition}
\Crefname{figure}{Figure}{Figures}
\Crefname{table}{Table}{Tables}
\Crefname{section}{Section}{Sections}
\Crefname{appendix}{Appendix}{Appendices}

\Crefname{problem}{Problem}{Problem}

\newcommand{\tr}{{{\mathsf T}}}
\newcommand{\Tr}{{{\mathbf{Tr}}}}

\newcommand{\LQR}{\mathtt{LQR}}
\newcommand{\Hinf}{\mathcal{H}_\infty}

\newcommand{\CM}{\mathcal{M}}
\newcommand{\CN}{\mathcal{N}}
\newcommand{\domuL}{\mathcal{L}_2^m[0,\infty)}
\newcommand{\domwL}{\mathcal{L}_2^p[0,\infty)}

\newcommand{\ABKstar}{A_{K^\star}}

\newcommand{\ABwLstar}{A_{L^\star}}
\newcommand{\bT}{{\mathbf{T}}}
\newcommand{\pol}{{\overline{p}}}

\newcommand{\pul}{{\underline{p}}}

\newcommand{\Kstar}{{K^\star}}
\newcommand{\Lstar}{{L^\star}}
\newcommand{\KstarTr}{(K^\star)^\mathsf{T}}
\newcommand{\LstarTr}{(L^\star)^\mathsf{T}}
\newcommand{\diag}{\mathrm{diag}}
\newcommand{\sdp}{\texttt{sdp}}

\begin{document}
\bstctlcite{IEEEexample:BSTcontrol}

\title{\bf Semidefinite Programming Duality in Infinite-Horizon Linear Quadratic Differential Games}

\author{Yuto Watanabe, 
% \IEEEmembership{Student member, IEEE}, 
Chih-Fan Pai, and Yang Zheng
% \IEEEmembership{Senior Member, IEEE}
\thanks{This work is supported by NSF ECCS-2154650, NSF CMMI 2320697, and NSF CAREER 2340713.}
\thanks{Y. Watanabe, C.-F. Pai,  and Y. Zheng are with the Department of Electrical
and Computer Engineering, University of California San Diego; \texttt{\{y1watanabe,cpai,zhengy\}@ucsd.edu}}}

\maketitle

\begin{abstract}
Semidefinite programs (SDPs) play a crucial role in control theory, traditionally as a computational tool. Beyond computation, the duality theory in convex optimization also provides valuable analytical insights and new proofs of classical results in control. In this work, we extend this analytical use of SDPs to study the infinite-horizon linear-quadratic (LQ) differential game in continuous time. 
Under standard assumptions, we introduce a new SDP-based primal-dual approach to establish the saddle point characterized by linear static policies in LQ games.
For this, we leverage the Gramian representation technique,
which elegantly transforms linear quadratic control problems into tractable convex programs.  
We also extend this duality-based proof to the $\Hinf$ suboptimal control problem. To our knowledge, 
this work provides the first primal-dual analysis using Gramian representations for the LQ game and $\Hinf$ control beyond LQ optimal control and $\mathcal{H}_\infty$ analysis. 
\end{abstract}

\section{Introduction}
\label{sec:introduction}
The linear quadratic (LQ) differential game is one fundamental problem in control \cite{bacsar2008h,engwerda2005lq}.
This problem deals with a non-cooperative dynamic decision-making process of two players, which is characterized by a linear dynamical system with a quadratic performance measure.
These two players are often called \textit{control} and \textit{disturbance}.
A pair of reasonable choices for both players is typically captured by the concept of \textit{Nash equilibria} (NE) or a \textit{saddle point}.  

Classical results in the LQ game \cite{bacsar2008h} include 
i) linear policies constitute an NE in the finite horizon case and ii) they may not form an NE in the infinite-horizon case unless a stronger assumption is imposed, but such a linear policy still works for the control player. 
A standard approach to establish these results is the \textit{completion-of-squares} technique \cite{mageirou1976values}. 
Thanks to these favorable properties, the LQ game and game-theoretic perspective have been extensively utilized for robust $\Hinf$ control which accounts for the worst-case disturbance  \cite{zhou1996robust,bacsar2008h,green2012linear}. 
Recently, the simple structure of the NE in LQ games also inspired new applications, such as a model-free NE-seeking approach called direct policy search \cite{zhang2019policy,hu2023toward,talebi2024policy}. 

It is known that semidefinite programming (SDP) plays an important role across a broad range of fundamental control problems, both theoretically and practically \cite{boyd1994linear}. 
In particular, SDPs offer implementable computations of foundational control problems, such as $\mathcal{H}_2$ and $\Hinf$ synthesis \cite{scherer1997multiobjective,gahinet1994linear}. These SDP characterizations further reveal benign nonconvex geometry in control problems \cite{zheng2023benign,zheng2024benign,tang2023analysis,hu2022connectivity}.  
Motivated by the fundamental importance of SDPs in control, several works \cite{balakrishnan2003semidefinite,gattami2015simple,you2015primal,bamieh2024linear,watanabe2025revisiting} have investigated their dual formulations and provided several new analytical insights. For example,  
several fundamental control problems, 
including the KYP lemma and linear quadratic regulator (LQR), 
were re-interpreted via a duality analysis in \cite{balakrishnan2003semidefinite}.
The authors of \cite{gattami2015simple,you2015primal}~introduced a primal-dual perspective for $\Hinf$ norm analysis and presented new computation methods.
Recently, the LQR was revisited via SDP duality in \cite{bamieh2024linear,watanabe2025revisiting}, which offers new analysis strategies. % were derived. 
However, for the LQ differential game and $\Hinf$ control,
the primal-dual perspective is still incomplete. 
A major difficulty is how to deal with min-max and max-min operations in LQ games, which makes~the existing results for $\Hinf$ analysis and LQR \cite{balakrishnan2003semidefinite,gattami2015simple,you2015primal,bamieh2024linear,watanabe2025revisiting} inapplicable.

In this work, we aim to provide a new primal-dual perspective for the infinite-horizon LQ differential game and $\Hinf$ control in continuous-time linear time-invariant (LTI) systems.
We present a new duality-based proof of the celebrated results in LQ games \cite{bacsar2008h}, and resolve the technical difficulty in
min-max and max-min procedures. 
With the standard assumption that an algebraic Riccati equation (ARE) has a stabilizing solution satisfying a regularity condition, we derive a saddle point, given by a pair of linear static policies. Our analysis clarifies several key ingredients in LQ games, such as the ARE and relation between the linear policies and Karush–Kuhn–Tucker (KKT) condition. 
A key technique in our proof is \textit{Gramian representation} \cite{skelton1997unified,you2015primal,watanabe2025revisiting}, by which we transform lower and upper bounds of the value of the LQ game into SDPs.
This technique also plays a pivotal role in the duality-based proofs of $\Hinf$ analysis \cite{gattami2015simple} and LQR \cite{bamieh2024linear,watanabe2025revisiting}.
We extend this Gramian representation technique to the LQ game involving min-max and max-min operations. Our SDP constructions for the LQ game reveal an interesting fact that  
the classical ARE appears in the dual SDPs, which allows us to establish the saddle point property. 
Finally, we also extend this strategy to $\Hinf$ control. 

To our knowledge, this is the first primal and dual SDP-based analysis for LQ games. This primal-dual approach further gives us a few insights: 
(i) By the Gramian representation, the infinite-dimensional variables (states, control, and disturbance) can be handled as finite-dimensional matrix variables corresponding to their energies. This transformation makes LQ games more accessible to standard convex optimization. 
(ii) The ARE solves the dual problems of two SDPs that are upper and lower bounds of the LQ game, which directly tells us the value of the game;
(iii) The explicit form of linear polices of the NE emerges from the \textit{complementary slackness} in the KKT condition~of~SDPs. 
We hope this primal and dual analysis will inspire further theoretical and practical developments in game-theoretic control. 

The remainder of this paper is organized as follows.
\cref{section:problem-formulation} presents the LQ game setup and  
the problem statement.
\cref{section:preliminaries} introduces SDP and Gramian-based analysis in control. % technical preliminaries.
In \cref{section:proof-main-theorem}, we present a new SDP duality-based proof, and we discuss $\Hinf$ control as an application in \cref{section:Hinf}. 
\cref{section:example} gives an illustrative~example.
Finally, \cref{section:conclusion} concludes this paper.

\paragraph*{Notations}
We use $\mathbb{S}^n_+$ to denote the set of positive semidefinite matrices. %are denoted by $\mathbb{S}^n_{++}$ ($\mathbb{S}^n_+$).
Given two symmetric matrices $M_1,M_2\in\mathbb{S}^n$, we use $M_1\prec (\preceq) M_2$ and $M_1\succ (\succeq) M_2$ when $M_1-M_2$ is negative (semi)definite and positive (semi)definite, respectively.
We use $\mathcal{L}_2^k[0,\infty)$ to denote the set of square-integrable (bounded energy) signals over time $[0,\infty)$ and with dimension $k$.
For a signal $u\in\mathcal{L}_2^k[0,\infty)$, we write its $2$-norm as
$\|u\|_2 = (\int_0^\infty \|u(t)\|^2dt)^{1/2}$.

\section{Preliminaries and Problem formulation}\label{section:problem-formulation}

We here review the zero-sum LQ game and a classical result characterizing the saddle-point property and present our problem statement of SDP duality in LQ games. 

\subsection{LQ differential game problem}
We consider a zero-sum game, where the system dynamics are characterized by a differential equation
\begin{align}\label{eq:dynamics}
        \dot{x}(t) = Ax(t) + Bu(t) + B_w w(t).
\end{align}
Here, the system state is $x(t)\in\mathbb{R}^n$, the inputs of~players 1 and 2 are $u(t) \in \mathbb{R}^m$ (control) and $w(t) \in \mathbb{R}^p$ (disturbance), respectively, 
and matrices $A\in\mathbb{R}^{n\times n}$, $B\in\mathbb{R}^{n\times m}$, $B_w\in\mathbb{R}^{n\times p}$ specify the system dynamics. 

To properly define the zero-sum game, we need to specify the information structures for both players. 
For the control player, we consider a \textit{closed-loop} information structure:
\begin{align*}
    u(t) = \mu(t;x(\tau),\,\tau\leq t),\quad t\geq 0
\end{align*}
where $\mu(t;x(\tau),\,\tau\leq t)$ is a Borel measurable function with piece-wise continuity for $t$ and Lipschitz continuity for $x$. 
It is allowed to depend on the entire state trajectory $x(\tau),\,\tau\leq t$.  
Analogously, we consider the same form of policies for the disturbance player $w$:
\begin{equation*}
    w(t) = \nu(t;x(\tau),\,\tau\leq t),\quad t\geq 0.
\end{equation*}
Let the policy space of the control player be  $\CM$, which~includes 
all such Borel measurable functions with piece-wise continuity for $t$ and Lipschitz continuity for $x$. Similarly, we use $\CN$ as the set of admissible policies for the disturbance player. 
For notational simplicity, 
we write $\mu(x(t))=Kx(t)$ and $ \nu(x(t))=Lx(t)$ for linear static state feedback policies.

We now specify a quadratic performance measure of the game $J_\gamma:\CM\times \CN\to\mathbb{R}$, defined as   
\begin{equation}\label{eq:J_gamma}
\begin{aligned}
    J_\gamma(\mu,\nu) %= \|z(t)\|_2^2-\gamma^2\|w(t)\|_2^2\\
    =%& 
    \int_0^\infty \left(
    x^\tr Qx+u^\tr Ru
    -\gamma^2\|w\|^2\right)dt,
\end{aligned}
\end{equation}
where $\gamma > 0$ is a fixed constant, $Q \succ 0$, and $R \succ 0$. 
The objective of Player 1 (resp., Player 2) is to minimize (resp., maximize) the function $J_\gamma(\cdot,\cdot)$ with respect to the function $\mu\in\CM$ (resp., $\nu\in\CN$).
Then, this constitutes a two-player non-cooperative zero-sum game.

Next,
we introduce the standard concept of \textit{a saddle point} (or \textit{Nash equilibrium (NE)}). 
The following $\pol^\star$ below is called the \textit{upper value} of the game:
\begin{subequations}\label{eq:LQ-game-value}
\begin{align}\label{eq:LQ-game-upper-value}
  \begin{aligned}
      \pol^\star=\inf_{\mu\in\CM
      }  \sup_{\nu\in\CN
      } \;
      J_\gamma(\mu,\nu)
      \quad
    \text{subject to} \;
    \cref{eq:dynamics},\, x(0)=x_0. 
  \end{aligned}
    \end{align}
On the other hand, the following $\pul^\star$ is called the \textit{lower value} of the game:
\begin{align}\label{eq:LQ-game-lower-value}
  \begin{aligned}
      \pul^\star=  \sup_{\nu\in\CN
      } 
      \inf_{\mu\in\CM
      }
      \;J_\gamma(\mu,\nu)
      \quad
    \text{subject to} \;
    \cref{eq:dynamics},\, x(0)=x_0.
  \end{aligned}
    \end{align}
\end{subequations}
In general, we have $\pul^\star \leq \pol^\star$. If we have 
\begin{equation}\label{eq:def-saddle}
    \pul^\star = \pol^\star = J_\gamma(\mu^\star,\nu^\star)
\end{equation}
with $(\mu^\star,\nu^\star)\in\CM\times \CN$, we say $p^\star = J_\gamma(\mu^\star,\nu^\star)$ is the value of the game, and the pair $(\mu^\star,\nu^\star)$ is a saddle point (or Nash equilibrium) of the game. 
We know that \cref{eq:def-saddle} is equivalent to
$$
J_\gamma(\mu^\star,\nu) \leq J_\gamma(\mu^\star,\nu^\star) \leq 
J_\gamma(\mu,\nu^\star), \;\; \forall \mu\in\CM, \nu\in\CN. 
$$

This paper addresses the following saddle-point problem. 

\vspace{2pt}

\begin{problem}\label{problem:1_saddle}
    Consider the infinite-horizon LQ differential game \cref{eq:LQ-game-value} with the performance measure $J_\gamma(\mu,\nu)$ in \cref{eq:J_gamma}. Find a saddle point (i.e., Nash Equilibrium) $(\mu^\star,\nu^\star) \in \CM\times \CN$.
\end{problem}
\vspace{2pt} 

This zero-sum LQ game and its variants have been widely studied; see the textbooks \cite{bacsar2008h,engwerda2005lq}. It also connects closely with $\mathcal{H}_\infty$ control \cite{zhou1996robust}. We here only emphasize the differences between open-loop and closed-loop information~structures. 

\begin{remark}[Open- vs closed-loop information
structures] 
%\yang{@Yuto, please add a short discussion here}
The information structure of the players' policy spaces~is a fundamental aspect of game-theoretic control.~To~capture non-cooperative interactions, these policy spaces must~be fully decoupled.
In 
this paper, we consider a \textit{closed-loop} information structure as $\CM$ and $\CN$ suitable for feedback policies. In contrast, 
feedforward policies~correspond~to open-loop information structures, where the policy depends only on time $t$ and the initial state $x_0$, but not~on~the actual state trajectory.
We refer to \cite{bacsar2008h} for other~variants, such as sampled data and delayed state structures.
\hfill $\square$
\end{remark}

\subsection{Algebraic Riccati equation and saddle point}

The algebraic Riccati equation (ARE) \cite{lancaster1995algebraic} plays an important role in characterizing the existence of a saddle point.
In particular, we consider the following~ARE
 \begin{gather}
     \label{eq:indefinite-ARE-gamma}
    \mathcal{R}(P) =0 \quad \text{with}\\
    \mathcal{R}(P) :=  A^\tr P + P A + Q
    - P BR^{-1}B^\tr P + \frac{1}{\gamma^2} P B_wB_w^\tr P.\nonumber
 \end{gather}

\vspace{-2pt}
We now make the following assumption.
\begin{assumption} \label{assumption:value-game}
    We consider $\gamma > 0$, and the ARE \cref{eq:indefinite-ARE-gamma} has a symmetric solution $P^\star$ such that a) is $P^\star$ positive~semidefinite and ensures $A-BR^{-1}B^\tr P^\star + \frac{1}{\gamma^2}B_wB_w^\tr P^\star$ is stable; b) we have $Q-\frac{1}{\gamma^2} P^\star B_wB_w^\tr P^\star \succ0$.
\end{assumption}

\vspace{2pt}
 
Then, we have the following theorem.   
\begin{theorem}
\label{theorem:main_theorem}
    Consider 
    the zero-sum LQ game \cref{eq:LQ-game-value} with an initial condition $x(0)\in\mathbb{R}^n$. With \cref{assumption:value-game}, 
   we have
       $
            \pul^\star = \pol^\star = x_0^\tr P^\star 
        $ 
    and the pair of linear static policies    \vspace{-2pt}    
\begin{subequations}\label{eq:u-w-pair}
\begin{align}\label{eq:u_optimal}
            \mu^\star (x(t)) =& -R^{-1}B^\tr P^\star x(t),\\
            \label{eq:w_optimal}
            \nu^\star(x(t))=& \frac{1}{\gamma^2}B_w^\tr P^\star x(t)
    \end{align}
\end{subequations}
    constitutes a saddle point of the LQ game \cref{eq:LQ-game-value}.
\end{theorem}

This theorem is classical; see \cite[Theorem 4.8]{bacsar2008h}.  
The existence of a stabilizing solution to the ARE \cref{eq:indefinite-ARE-gamma} in \Cref{assumption:value-game} is standard, which guarantees the upper and lower values
are the same, i.e., $\pul^\star = \pol^\star = x_0^\tr P^\star x_0$.  
The condition $Q-\frac{1}{\gamma^2} P^\star B_wB_w^\tr P^\star \succ0$ is also needed to establish the saddle point property with the closed-loop information structure; otherwise, the pair
\cref{eq:u-w-pair} may not be a saddle point;~see \cite[Example 4.1]{bacsar2008h} and \cite{mageirou1976values}.
A similar assumption can be found in \cite{zhang2019policy}. 
In our later analysis, the assumption of $Q-\frac{1}{\gamma^2} P^\star B_wB_w^\tr P^\star \succ0$ is redundant for analyzing the upper value $\pol^\star$, as expected from \cite[Theorem 4.8]{bacsar2008h}. 
We will present a numerical illustration on the role of $Q-\frac{1}{\gamma^2} P^\star B_wB_w^\tr P^\star \succ0$ in \cref{section:example}.

\subsection{Problem statement: SDP duality in LQ games}

Under \cref{assumption:value-game}, a standard proof to establish \cref{theorem:main_theorem} is the \textit{completion-of-squares} \cite{mageirou1976values,engwerda2005lq,bacsar2008h,green2012linear}.
A key step in this approach is to add and subtract
$
    \int_0^\infty x^\tr(t) P^\star x(t) dt
$ 
for $J_\gamma$ in \cref{eq:J_gamma}. After some calculations, $J_\gamma$ can be rewritten~as
\begin{align*}
&J_\gamma(\mu,\nu)
=
x_0^\tr P^\star  x_0
-\int_0^\infty \gamma^2\left\|w(t)-\nu^\star\left(x(t)\right)\right\|^2 dt
    \\
    &
    +\int_0^\infty
    \left(u(t)-\mu^\star\left(x(t)\right)\right)^\tr R\left(u(t)-\mu^\star\left(x(t)\right)\right) 
    dt.
\end{align*}
We can see that $x_0^\tr P^\star  x_0$ is the value of the game, and \cref{eq:u-w-pair} is a saddle point. Although this strategy is simple, it is not very useful for further analysis (e.g., the existence of the stabilizing solution),
and such analysis requires more~involved arguments, e.g., taking the limit of the finite-horizon case~\cite{bacsar2008h}. 

In this work,
we aim to 
reveal the role of SDPs and duality-based analysis for LQ games.
In particular, we 
establish \cref{theorem:main_theorem} and solve $\Hinf$ control using SDP duality. 
For LQR and $\Hinf$ analysis,
a primal-dual perspective has already appeared in \cite{bamieh2024linear,watanabe2025revisiting,gattami2015simple,you2015primal}, which offers extra theoretical insights and new computation methods (e.g., a new proof for the LQR optimal gain \cite{bamieh2024linear,watanabe2025revisiting}, structured $\Hinf$ analysis \cite{gattami2015simple,you2015primal}).
One notable benefit of this perspective is the ease of analyzing the relation between control-theoretic regularity conditions (e.g., controllability and observability) and the existence of optimal/feasible solutions,
using well-established duality results.
However, such a primal and dual perspective for the LQ game and $\Hinf$ control remains open. 
One challenge lies in handling the min-max and max-min problems in \cref{eq:LQ-game-value}, which represents perhaps the most significant gap between LQR and LQ games. 
With this motivation, our main objectives are twofold:
\begin{enumerate}
    \item By leveraging SDPs and duality,
    we provide an alternative and self-contained proof for \cref{theorem:main_theorem}. % and reconcile the gap; 
    \item We resolve the challenge and clarify the role of the key building blocks, such as the ARE \cref{eq:indefinite-ARE-gamma} and static linear policies \cref{eq:u-w-pair}, from the perspective of SDP duality. 
\end{enumerate}

\section{SDP Duality and Gramian Representation in Control}\label{section:preliminaries}

In this section, we first review standard duality results in SDPs. Then, we introduce the notion of 
\textit{Gramian} representation\footnote{Its stochastic variant is called \textit{covariance} representation.} that captures all quadratic information of the trajectories from an LTI system \cite{skelton1997unified}.

\subsection{Duality in SDPs}

The standard primal SDP is a problem of the form:
\begin{align}\label{eq:sdp-primal-standard}
    \!\!
    \begin{aligned}
p^*=\min_{X} &\; \langle C,X\rangle \quad
\text {subject to} \quad  
\mathcal{A}(X)=b,\,X \in\mathbb{S}^r_+,
\end{aligned}
\end{align}
%\end{subequations}
where $b\in\mathbb{R}^m,\,C\in\mathbb{S}^r$ and
$\mathcal{A}:\mathbb{S}^r\to\mathbb{R}^m$ is a linear map.
% ,b_i\in\mathbb{R},\,i=1,\ldots,m$.
Its Lagrange dual problem is 
\begin{align}\label{eq:sdp-dual-standard}
    \begin{aligned}
d^*=
\max_{y\in\mathbb{R}^m}& \quad b^\tr y\quad
 \text {subject to } \quad 
 C - \mathcal{A}^*(y) \in\mathbb{S}^r_+, 
\end{aligned}
\end{align}
where $\mathcal{A}^*:\mathbb{R}^m\to\mathbb{S}^r$ is the adjoint of $\mathcal{A}(\cdot)$ defined by
$
 \left\langle
    y,\mathcal{A}(X)
    \right\rangle=
    \left\langle
    \mathcal{A}^*(y),X
    \right\rangle, \forall y \in \mathbb{R}^m, X \in \mathbb{S}^r.  
$ 

The KKT optimality condition for this pair of primal-dual SDPs \cref{eq:sdp-primal-standard,eq:sdp-dual-standard} consists of three properties:
\begin{itemize}
    \item Primal feasibility: 
    $X \in \mathbb{S}_{+}^r$ and
    $\mathcal{A}(X)=b$; 
    \item Dual feasibility: $C - \mathcal{A}^*(y)\in\mathbb{S}^r_+$;
    \item Complementary slackness: $\left\langle 
    X
    ,C-\mathcal{A}^*(y)\right\rangle=0$.
\end{itemize}

We next present a version of the duality in SDPs, and 
%This lemma is our main tool for establishing \cref{theorem:main_theorem}.
we refer to \cite{vandenberghe1996semidefinite} for further discussions on strong duality. 

\begin{lemma}\label{theorem:SDP-duality}
    For the pair of primal and dual SDPs \cref{eq:sdp-primal-standard,eq:sdp-dual-standard},
    the following statements hold.
    \begin{enumerate}
        \item Weak duality: $d^*\leq p^*$.
    \item Strong duality:
    If there exists a primal-dual feasible solution pair $\left(X^\star,y^\star\right)$ satisfying the complementary slackness
    $\left\langle 
    X^\star
    ,C-\mathcal{A}^*(y^\star)\right\rangle=0$,
    then strong duality holds, i.e., $d^*= p^*$, and
    $X^\star$ and $y^\star$ are optimal for the primal and dual problems, respectively.
    \end{enumerate}
\end{lemma}

Weak duality always holds for any primal and dual problems, not limited to SDPs.
The second statement with the KKT condition directly verifies the optimality of a pair of primal-dual solutions  $\left(X,y\right)$. Under mild conditions (such as strict feasibility and finite optimal values), the KKT condition and strong duality for SDPs are always guaranteed \cite{vandenberghe1996semidefinite}. 

\subsection{Gramian representation in control problems}\label{subsection:gramian-LQR}

Here, we introduce the notion of Gramian representation, which plays a fundamental role in SDP-based analysis for control problems. 
For simplicity, consider an LTI system
\begin{equation}\label{eq:LTI-Gramian}
    \dot{x}(t) = Ax(t) + Bu(t),\quad x(0)=x_0,
\end{equation}
where $x(t) \in \mathbb{R}^n$ and $u(t) \in \mathbb{R}^m$. 
For any input $u(t)\in \mathcal{L}_2^m[0,\infty)$ that ensures $x(t)\in\mathcal{L}_2^n[0,\infty)$, we can define a positive semidefinite matrix  
\begin{equation} \label{eq:gramian-control}
    Z
=\displaystyle\int_0^\infty 
    \begin{bmatrix}x(t)\\u(t)\end{bmatrix}
    \begin{bmatrix}x(t)\\u(t)\end{bmatrix}^\tr dt
    \in \mathbb{S}_+^{n+m},  
\end{equation}
which encodes the quadratic information of the trajectory $x(t),u(t)$. This matrix \cref{eq:gramian-control} is also known as the Gramian matrix. Accordingly, we define the following set 
\begin{align} \label{eq:setV}
    \begin{aligned}
   \mathcal{V} \!= \! 
   \bigcup_{u\in\domuL}
   \left\{
    Z = \cref{eq:gramian-control}
    \in \mathbb{S}^{n+m} 
    \middle    |
x(t)\in\mathcal{L}_2^n[0,\infty) \text{ for \cref{eq:LTI-Gramian}}
    \right\}.
\end{aligned} 
\end{align}
\normalsize
This set contains the full quadratic information of all the stable trajectories starting from $x(0)=x_0$ \cite{skelton1997unified}.
It is then possible to express a linear quadratic control problem by replacing $x(t)$ and $u(t)$ by $Z\in\mathcal{V}$.
For example, consider the LQR problem, i.e., minimizing the quadratic cost:
\begin{equation}\label{eq:LQR}
\begin{aligned}
 p_\LQR^\star\!= \!  \min_{u\in\domuL}\
\int_0^\infty \!
\left(x(t)^\tr Qx(t)
\!+\! u(t)^\tr Ru(t)
\right)dt
\end{aligned}
\end{equation}
subject to \cref{eq:LTI-Gramian}.
If $R\succ 0$ and $(Q^{1/2},A)$ is detectable,
it follows from \cite[Lemma 14.1]{zhou1996robust} that \cref{eq:LQR} can 
 be rewritten as
\begin{equation}\label{eq:LQR-gram}
   p_\LQR^\star:= \min_{Z\in\mathcal{V}}
 \left\langle Z,\diag(Q,R) \right\rangle.
\end{equation}
Notice that
this reformulation replaces the quadratic cost in \cref{eq:LQR} by the linear cost \cref{eq:LQR-gram} with respect to $Z$.

However, \cref{eq:LQR-gram} is not immediately tractable in its current form since $\mathcal{V}$ lacks an explicit characterization.
%To handle the set and \cref{eq:LQR-gram} 
For duality-based analysis, we introduce an outer convex approximation 
\begin{equation}  \label{eq:setV-outer}
\mathcal{V}^\sdp
    = \left\{
    Z \in \mathbb{S}^{n+m} \middle| \begin{aligned}
    &x_0x_0^\tr + AZ_{11}+BZ_{12}^\tr + \\
    &\qquad   Z_{11}A^\tr + Z_{12}B^\tr  = 0, \\
    &\quad Z = 
    \begin{bmatrix}Z_{11} & Z_{12}\\ Z_{12}^\tr & Z_{22}\end{bmatrix} \succeq 0
        \end{aligned}
    \right\}.
\end{equation}
It is known that $\mathcal{V} \subseteq \mathcal{V}^\sdp$ \cite[Lemma 10]{watanabe2025revisiting}. 
We then have
\begin{equation}\label{eq:LQR-gram_SDP}
\begin{aligned}
    p_\LQR^\star\geq 
   p_\LQR^\sdp
:= &\min_{Z\in\mathcal{V}^\sdp}\quad
 \left\langle Z,\diag(Q,R) \right\rangle, 
\end{aligned}
\end{equation}
which is a convex SDP of the form \cref{eq:sdp-primal-standard}.
Thus, the Gramian representation \cref{eq:gramian-control} allows us to handle the LQR problem \cref{eq:LQR} in the SDP framework. We 
can further derive the optimal solution to LQR \cref{eq:LQR} using \cref{eq:LQR-gram} and \cref{eq:LQR-gram_SDP} and their duals~\cite{watanabe2025revisiting}.

It is standard to derive the Lagrange dual problem of \cref{eq:LQR-gram_SDP}:  
\begin{equation}\label{eq:dual-LQR}
\begin{aligned}
p_\LQR^\sdp
\geq d_\LQR =&\max_{P \in \mathbb{S}^n}\quad
x_0^\tr P x_0
\\
\text { subject to }\quad&
\begin{bmatrix}
A^{\top} P+P A+Q & P B \\
B^{\top} P & R    
\end{bmatrix}
\succeq 0 .
\end{aligned}
\end{equation}
By a celebrated comparison theorem \cite{willems1971least},
an optimal solution to \cref{eq:dual-LQR} is given by
the unique stabilizing solution $P=\hat{P}^\star\succeq 0$ to the ARE below
\begin{equation}\label{eq:ARE-LQR}
A^\tr P + P A + Q
    - P BR^{-1}B^\tr P 
          =0. 
\end{equation}
Such a stabilizing solution is guaranteed to exist if $(A,B)$ is stabilizable and $(Q^{1/2},A)$ is detectable.

Now,
using 
a stable input
$u(t)=\hat{K}^\star x(t)$ with
$\hat{K}^\star = - R^{-1}B^\tr \hat{P}^\star$, % to the original LQR \cref{eq:LQR}, 
we construct a feasible solution $\hat{Z}^\star$ to \cref{eq:LQR-gram_SDP}: 
\begin{align}\label{eq:Zstar-LQR}
\begin{aligned}
\hat{Z}^\star &=\int_0^\infty \begin{bmatrix}
    x(t)\\
    \hat{K}^\star x(t)
\end{bmatrix}
\begin{bmatrix}
    x(t)\\
    \hat{K}^\star x(t)
\end{bmatrix}^\tr dt
= \begin{bmatrix}
    I\\
    \hat{K}^\star 
\end{bmatrix}
Z_{11}^\star
\begin{bmatrix}
    I\\
    \hat{K}^\star 
\end{bmatrix}^\tr
\end{aligned}
\end{align}
with $Z_{11}^\star:=\int_0^\infty x(t)x(t)^\tr dt$.
Then, we can verify the complementary slackness in \Cref{theorem:SDP-duality} as
\begin{align*}
&\left\langle 
\hat{Z}^\star,\left[\begin{array}{cc}
A^{\top} \hat{P}^\star+\hat{P}^\star A+Q & \hat{P}^\star B \\
B^{\top} \hat{P}^\star & R
\end{array}\right]
\right\rangle\\
\overset{\cref{eq:ARE-LQR}}{=}&
\left\langle 
\hat{Z}^\star,\left[\begin{array}{cc}
\hat{P}^\star BR^{-1}B^\tr \hat{P}^\star & \hat{P}^\star B \\
B^{\top} \hat{P}^\star & R
\end{array}\right]
\right\rangle\\
=&
\left\langle 
\hat{Z}^\star,
\begin{bmatrix}
(\hat{K}^\star)^\tr R^{1/2}\\
-R^{1/2}
\end{bmatrix}
\begin{bmatrix}
(\hat{K}^\star)^\tr R^{1/2}\\
-R^{1/2}
\end{bmatrix}^\tr
\right\rangle \overset{\cref{eq:Zstar-LQR}}{=}0.
\end{align*}
Thus, \Cref{theorem:SDP-duality} implies that $\hat{Z}^\star$ is an optimal solution to \cref{eq:LQR-gram_SDP}.
Since \cref{eq:LQR-gram_SDP} gives a lower bound of $p_\LQR^\star$ in the 
LQR \cref{eq:LQR}, we conclude that
$u(t)=\hat{K}^\star x(t) = - R^{-1}B^\tr \hat{P}^\star x(t)$ is globally optimal to \cref{eq:LQR}.
Further details on strong duality in LQR can be found in \cite{watanabe2025revisiting}; also see \cite{bamieh2024linear,gattami2015simple,you2015primal}.  

The SDP-duality analysis offers an alternative way to solve the classical LQR. 
We note that the dual \cref{eq:dual-LQR} is solved by the ARE \cref{eq:ARE-LQR}, and the form of optimal gain $\hat{K}^\star$ is determined by the KKT condition, especially complementary slackness.
We will extend this SDP duality approach to analyze the upper and lower values $\pol^\star$ and $\pul^\star$ in \cref{section:proof-main-theorem}.

\section{Nash Equilibrium via SDP duality}\label{section:proof-main-theorem}

We here present a new proof of \cref{theorem:main_theorem} based on the SDP duality and the Gramian technique. 
To prove the saddle point property, we take the following steps:
\begin{enumerate}
\item[a)] Via the Gramian representation,
we derive two SDPs offering upper and lower bounds of $\pol^\star$ and $\pul^\star$, resp.;
\item[b)] We show that the ARE \cref{eq:indefinite-ARE-gamma} appears in the duals of both SDPs in step (a). Especially, $P^\star$ in \cref{eq:indefinite-ARE-gamma} is feasible to their duals, and this confirms $\pol^\star=\pul^\star = x_0^\tr P^\star x_0$;
\item[c)] By KKT analysis in the primal and dual SDPs, we establish
$\pol^\star = \max_{\nu\in\CN} J_\gamma(\mu^\star,\nu)$ with
\begin{subequations}
    \begin{align}\label{eq:argmax=wstar}
            \nu^\star\in
        {\arg \max}_{\nu\in\CN} J_\gamma(\mu^\star,\nu),
        \end{align}
and
$\pul^\star = \min_{\mu\in\CM} J_\gamma(\mu,\nu^\star)$ 
    with
        \begin{align}\label{eq:argmin=ustar}
            \mu^\star\in
        {\arg \min}_{\mu\in\CM} J_\gamma(\mu,\nu^\star),
        \end{align}
\end{subequations}
where $\mu^\star,\nu^\star$ are the linear static polices in \cref{eq:u-w-pair}. 
\end{enumerate}
Throughout the proof, we use the standard duality result in \cref{theorem:SDP-duality}.
Similar to the LQR case in \cref{subsection:gramian-LQR}, the KKT conditions play the central role in step (c), and the linear static policies \cref{eq:u-w-pair} appear in solving the KKT with the ARE \cref{eq:indefinite-ARE-gamma}. 
In the following, we denote \begin{equation*}
    (K^\star,L^\star) = \left(-R^{-1}B^\tr P^\star,\frac{1}{\gamma^2}B_w^\tr P^\star
    \right).
\end{equation*}

\subsection{Step (a): SDPs as upper and lower bounds} \label{subsection:SDP-upper-lower-values}
We first derive an upper bound of $\pol^\star$ in \cref{eq:LQ-game-upper-value} using the Gramian representation. 
Restricting the control player's policy as $\mu(x(t))=\Kstar x(t)$ leads to following upper bound: \vspace{-10pt}
\begin{subequations} \label{eq:LQ-game_Upper_Kstar}
\begin{align}
       \pol^\star \leq  \pol_{\Kstar}=\sup_{\nu\in\CN} \quad& 
       J_\gamma(\Kstar x,\nu)
       \\
    \text{subject to} \quad&  \dot x(t) = \ABKstar x(t)+B_ww(t),\, \label{eq:LQ-game_Upper_Kstar-b}\\
    &x(0)=x_0,\,
\end{align}
\end{subequations}
where $\ABKstar:=A+B\Kstar$ and $w(t) = \nu(t;x(\tau),\,\tau\leq t).$
This problem \cref{eq:LQ-game_Upper_Kstar} resembles the structure in the LQR \cref{eq:LQR},~except that the disturbance player aims to maximize the cost in \cref{eq:LQ-game_Upper_Kstar}.

We next drive an SDP for \cref{eq:LQ-game_Upper_Kstar} to further upper bound $\pol_{\Kstar}$. 
It is known that $\Kstar$ stabilizes the LTI system \cite[Theorem 6.3.1]{green2012linear}, i.e., $\ABKstar$ is stable. The existence of a stabilizing solution $P^\star$ also guarantees the boundedness of $\pol^\star$ and $\pul^\star$ \cite[Theorem 4.8 (iii)]{bacsar2008h}. Accordingly, we can consider $w$ to have finite energy without loss of generality. This also ensures $x\in\mathcal{L}_2^n[0,\infty)$ due to the stability of $\ABKstar$. Then, for any $w$ with finite energy, the  Gramian matrix below is well-defined 
\begin{equation} \label{eq:gramian-disturbance}
Z = \begin{bmatrix}
    Z_{11} &Z_{12}\\
    Z_{12}^\tr &Z_{22}
\end{bmatrix}
= \int_0^\infty \begin{bmatrix}
    x(t)\\
    w(t)
\end{bmatrix}
\begin{bmatrix}
    x(t)\\
    w(t)
\end{bmatrix}^\tr dt \in \mathbb{S}^{n + p}_+. 
\end{equation}
We then define the following set 
\begin{align*}
   \mathcal{V}_{w} \!= \! \left\{
    Z = \cref{eq:gramian-disturbance}
    \in \mathbb{S}^{n+p} 
    \middle    |
    \begin{aligned} 
    w(t)\!\in\!\domwL \text{ in \cref{eq:LQ-game_Upper_Kstar-b}}
    \end{aligned}
    \right\}, 
\end{align*}
which allows us to equivalently reformulate \cref{eq:LQ-game_Upper_Kstar} as  
\begin{equation*}
    \pol_\Kstar= \sup_{Z\in\mathcal{V}_{w}}
 \left\langle Z,\diag(Q^\star,-\gamma^2I) \right\rangle,
\end{equation*}
where $Q^\star = Q+\KstarTr R\Kstar$.

The set $\mathcal{V}_{w}$ is not tractable. Similar to the outer approximation in \cref{eq:setV-outer}, we derive the following upper bound  
\begin{align}\label{eq:LQ-game_Upper_Kstar-gram}
\begin{aligned}
\pol_\Kstar
\leq \pol_{\Kstar}^\sdp:=&\sup_{Z\in \mathbb{S}^{n+p}_+}\quad
\left\langle Z,\diag(Q^\star,-\gamma^2I) \right\rangle\\
\text{subject to}& \quad x_0x_0^\tr +  \ABKstar Z_{11}+ B_wZ_{12}^\tr + \\
&\qquad \quad (\ABKstar Z_{11}+ B_wZ_{12}^\tr)^\tr  
=0,
\end{aligned}
\end{align}
which is an SDP. The Gramian representation \cref{eq:gramian-disturbance} enables the upper bound $\pol_{\Kstar}^\sdp \geq \pol^\star$ by solving the SDP \cref{eq:LQ-game_Upper_Kstar-gram}. 

\vspace{2pt}
\begin{proposition}  \label{proposition:upper-value}
    Consider 
    the LQ game \cref{eq:LQ-game-value} with an initial condition $x(0)\in\mathbb{R}^n$. With \cref{assumption:value-game}, the SDP \cref{eq:LQ-game_Upper_Kstar-gram} returns an upper bound for the upper value, i.e., $\pol^\star \leq \pol_{\Kstar}^\sdp $. 
\end{proposition}

\begin{proof}
    This result directly follows the discussions above by combining \cref{eq:LQ-game_Upper_Kstar} with \cref{eq:LQ-game_Upper_Kstar-gram}. 
\end{proof}

\vspace{2pt}

Similarly, we can construct a lower bound for the~lower value $\pul^\star$ in \cref{eq:LQ-game-lower-value}.
In particular, restricting the disturbance play's policy $\nu(x(t)) = L^\star x(t)\in\CN$ leads~to % the following problem:
% gives a lower bound of $\pul^\star$:
\begin{subequations}\label{eq:LQ-game_Lower-L}
\begin{align}
\pul^\star \geq
\pul_\Lstar=
\inf_{\mu\in\CM} \quad & J_\gamma(\mu,\Lstar x)\\
\text{subject to}\quad&
\dot{x}(t) = \ABwLstar x(t)+Bu(t), \label{eq:LQ-game_Lower-L-b}\\
&x(0)=x_0,
\end{align}
\end{subequations}
where $\ABwLstar=A+B_w\Lstar$ and $u(t)=\mu(t;x(\tau),\,\tau\leq t)$. This problem \cref{eq:LQ-game_Lower-L} is in the same form as the LQR \cref{eq:LQR}. Moreover, since $Q-\frac{1}{\gamma^2}P^\star B_wB_w^\tr P^\star\succ 0$,
the value $\pul_\Lstar$ in \cref{eq:LQ-game_Lower-L} is finite if and only if both $x$ and $u$ have finite energy. We thus consider only the case where $x\in\mathcal{L}_2^n[0,\infty)$ and $u\in\domuL$.
Then, we can introduce the set of Gramian matrices as
$$
    \begin{aligned}
   \mathcal{V}_u \!\!= 
   \bigcup_{u\in\domuL}
   \!\! \left\{
    Z = \cref{eq:gramian-control}
    \in \mathbb{S}^{n+m} 
    \middle    |
    \begin{aligned} 
    &x(t)\in\mathcal{L}_2^n[0,\infty) \\
    &\text{ in \cref{eq:LQ-game_Lower-L-b}}
    \end{aligned}
    \right\}, 
\end{aligned} 
$$
which
provides the equivalent problem as
\begin{equation*}
\pul_\Lstar=
\inf_{Z\in\mathcal{V}_u} 
 \left\langle Z,\diag(Q-\gamma^2 \LstarTr \Lstar,R) \right\rangle.
\end{equation*}
Similar to the outer approximation in \cref{eq:setV-outer}, we derive the following lower bound 
\begin{align}
        \underline{p}_\Lstar\geq \underline{p}_\Lstar^\sdp=& \inf_{Z\in\mathbb{S}_+^{n+m}} \quad \left\langle
        \begin{bmatrix}
            Q -\gamma^2 \LstarTr \Lstar & 0\\
            0 & R
        \end{bmatrix},Z\right\rangle \nonumber \\
    \text{subject to}\quad& x_0x_0^\tr +  \ABKstar Z_{11}+ B_wZ_{12}^\tr +  \label{eq:LQ-game_Lower_L-gram}\\
&\qquad \quad (\ABKstar Z_{11}+ B_wZ_{12}^\tr)^\tr  \nonumber 
=0.
    \end{align}
To summarize, we have obtained the following result. 

\vspace{2pt}
\begin{proposition} \label{proposition:lower-value}
    Consider 
    the LQ game \cref{eq:LQ-game-value} with an initial condition $x(0)\in\mathbb{R}^n$. With \cref{assumption:value-game}, the SDP \cref{eq:LQ-game_Lower_L-gram} returns a lower bound for the lower value, i.e., $\underline{p}_\Lstar^\sdp \leq \pul^\star$. 
\end{proposition}

\vspace{2pt}

After fixing the control or disturbance policy, the Gramian representation technique allows us to construct two SDPs \cref{eq:LQ-game_Lower_L-gram,eq:LQ-game_Upper_Kstar-gram} to bound the upper and lower values of the LQ game \cref{eq:LQ-game-value}.  
It turns out that the duals of \cref{eq:LQ-game_Lower_L-gram,eq:LQ-game_Upper_Kstar-gram} are related by the same ARE \cref{eq:indefinite-ARE-gamma}. This fact further establishes that $\underline{p}_\Lstar = \pol_{\Kstar}$. We provide the details in the next section.

\subsection{Step (b): Equal upper and lower values by SDP duality} \label{subsection:duality-analysis}

We establish the following result in this section.
\begin{proposition} \label{propostion:SDPs}
    Consider the SDPs \cref{eq:LQ-game_Lower_L-gram,eq:LQ-game_Upper_Kstar-gram}. With  \cref{assumption:value-game}, we have $\underline{p}_\Lstar^\sdp \geq x_0^\tr P^\star x_0 \geq \pol_{\Kstar}^\sdp$. 
\end{proposition}

\vspace{2pt}

From \Cref{proposition:upper-value,proposition:lower-value}, we have already established 
$$
\underline{p}_\Lstar^\sdp \leq \pul^\star \leq \pol^\star \leq  \pol_{\Kstar}^\sdp.
$$
Thus, together with \cref{propostion:SDPs}, the upper and lower values of the game must be the same, i.e., $\pul^\star = \pol^\star = x_0^\tr P^\star x_0$. 

We prove \cref{propostion:SDPs} by investigating the duals of \cref{eq:LQ-game_Lower_L-gram,eq:LQ-game_Upper_Kstar-gram}, both of which are related by the same ARE \cref{eq:indefinite-ARE-gamma}. 
We derive the dual problem of \cref{eq:LQ-game_Upper_Kstar-gram} as
\begin{align}\label{eq:LQ-game_Upper_Kstar-gram-dual}
        \begin{aligned}
\overline{d}_\Kstar=\inf_{P\in\mathbb{S}^n} &\; 
           x_0^\tr P  x_0
            \\
    \text{subject to}&\;
    \begin{bmatrix}
        \ABKstar^\tr P+P\ABKstar 
        + Q^\star
        & PB_w\\
        B_w^\tr P & - \gamma^2I
    \end{bmatrix}
    \preceq 0,
\end{aligned}
    \end{align}
and we have $\pol_\Kstar^\sdp\leq\overline{d}_\Kstar$ by weak duality.
The derivation is almost the same as 
\cite[Chapter 4.1.1]{you2016direct} and
\cite[Lemma 1]{watanabe2025revisiting}.
The constraint in \cref{eq:LQ-game_Upper_Kstar-gram-dual} is closely related to an \textit{algebraic Riccati inequality} (ARI) \cite{willems1971least}.
    Applying Schur complement, the constraint in \cref{eq:LQ-game_Upper_Kstar-gram-dual} is equivalent to 
    \begin{equation*}
        \ABKstar^\tr P+P\ABKstar 
        + Q^\star
        + \frac{1}{\gamma^2}PB_wB_w^\tr P \preceq  0.
    \end{equation*}
    From $Q^\star =Q+\KstarTr R\Kstar$,
 completing the squared term of $K^\star$ allows us to rewrite this inequality as
    \begin{equation}
     \label{eq:primal-constraint-completed}
    \begin{aligned}
           \mathcal{R}(P)\preceq 
           -(K^\star+R^{-1}B^\tr P)^\tr R(K^\star+R^{-1}B^\tr P).
    \end{aligned}
\end{equation}
It is clear that
    the solution $P^\star$ from
    the ARE $\mathcal{R}(P^\star)=0$ in \cref{eq:indefinite-ARE-gamma}    
    is feasible to \cref{eq:primal-constraint-completed} and thus also feasible to the dual SDP \cref{eq:LQ-game_Upper_Kstar-gram-dual}. This~implies
\begin{equation}\label{eq:value_inequality_upper_side}
    \pol_\Kstar^\sdp \leq \overline{d}_\Kstar
    \leq x_0^\tr P^\star x_0.
    \end{equation}

We next drive the dual problem of \cref{eq:LQ-game_Lower_L-gram} as
    \begin{align} \label{eq:LQ-game_Lower_L-gram-dual}
    \underline{d}_\Lstar=
    \sup_{P\in\mathbb{S}^n} &\;\;  x_0^\tr P x_0\\
        \text{subject to}&\; 
        \begin{bmatrix}
    \ABwLstar^\tr P\! + \!P\ABwLstar \!+ \!Q-\gamma^2 \LstarTr \Lstar& PB \\ B^\tr P & R
    \end{bmatrix}\succeq 0. \nonumber
    \end{align}
Since $R \succ 0$, the Schur complement shows that the constraint in \cref{eq:LQ-game_Lower_L-gram-dual} is equivalent to 
    \begin{align*}
       \ABwLstar^\tr P + P\ABwLstar + Q-\gamma^2 \LstarTr \Lstar
       -PBR^{-1}B^\tr P \succeq 0,
    \end{align*}
which can be rewritten as
\begin{equation} \label{eq:dual-constraint-completed}
\begin{aligned}
        \mathcal{R}(P)\succeq
       \gamma^2
       \left(
       L^\star-\frac{1}{\gamma^2}B_w^\tr P \right)^\tr
      \left(
       L^\star -\frac{1}{\gamma^2}B_w^\tr P \right).
    \end{aligned}
\end{equation}
    Here, we have completed the squared term of $L^\star$.
Now, we observe that the solution $P^\star$ to the ARE $\mathcal{R}(P^\star)=0$ in \cref{eq:indefinite-ARE-gamma} also satisfies this inequality, implying that $P^\star$ is feasible to the dual SDP \cref{eq:LQ-game_Lower_L-gram-dual}.
Thus, we obtain
\begin{equation}\label{eq:value_inequality_lower_side}
\pul_\Lstar^\sdp
\geq \underline{d}_\Lstar \geq x_0^\tr P^\star x_0 .
\end{equation}

Finally, combining \cref{eq:value_inequality_upper_side} and \cref{eq:value_inequality_lower_side} leads to \cref{propostion:SDPs} as
$$
\underline{p}_\Lstar^\sdp \geq x_0^\tr P^\star x_0 \geq \pol_{\Kstar}^\sdp . 
$$

\begin{remark}[Role of the ARE in LQ control]
    It is very interesting that the dual SDPs \cref{eq:LQ-game_Upper_Kstar-gram-dual,eq:LQ-game_Lower_L-gram-dual} are related by the same ARE \cref{eq:indefinite-ARE-gamma}, as shown in \cref{eq:primal-constraint-completed,eq:dual-constraint-completed}. Essentially, the ARE  \cref{eq:indefinite-ARE-gamma} characterizes the optimal solution to both dual SDPs \cref{eq:LQ-game_Upper_Kstar-gram-dual,eq:LQ-game_Lower_L-gram-dual}, which is the key to establishing \cref{propostion:SDPs}. We note that the ARE \cref{eq:ARE-LQR} also solves the dual SDP arising from the standard LQR \cref{eq:dual-LQR}. \hfill $\square$  
\end{remark}

\subsection{Step (c): Establishing saddle point \cref{eq:u-w-pair} via KKT analysis}

We finally establish the saddle point \cref{eq:u-w-pair} via KKT analysis for the primal and dual SDPs in \cref{subsection:duality-analysis}. 
We first establish
$\pol^\star = \max_{\nu\in\CN} J_\gamma(\mu^\star,\nu)$
and show \cref{eq:argmax=wstar}, i.e., $\nu^\star(x(t))=L^\star x(t)$ is a maximizer. 
To achieve this, we use the KKT condition for primal and dual SDPs
\cref{eq:LQ-game_Upper_Kstar-gram,eq:LQ-game_Upper_Kstar-gram-dual}: \vspace{-12pt}
\begin{subequations}\label{eq:KKT-upper}
\begin{align}
x_0x_0^\tr
+\ABKstar Z_{11}+ B_wZ_{12}^\tr \qquad \qquad  \qquad & \qquad \nonumber \\ 
+\left(\ABKstar Z_{11}+ B_wZ_{12}^\tr\right)^\tr
=0,
Z&\in\mathbb{S}_+^{n+p}, \label{eq:KKT-upper-primal-feasibility}
\\
\label{eq:KKT-upper-dual-feasibility}
\begin{bmatrix}
        \ABKstar^\tr P+P\ABKstar 
        + Q^\star
        % +Q+\KstarTr R\Kstar
        & PB_w\\
        B_w^\tr P & - \gamma^2I
    \end{bmatrix}&\succeq 0,\\
\left\langle 
Z,\begin{bmatrix}
\ABKstar^\tr P+P\ABKstar 
        + Q^\star
        % +Q+\KstarTr R\Kstar
        & PB_w\\
        B_w^\tr P & - \gamma^2I
    \end{bmatrix}
\right\rangle&=0.
\label{eq:KKT-upper-complementary-slackness}
\end{align}
\end{subequations}
We will construct a primal-dual solution pair $(Z,P)$ satisfying \cref{eq:KKT-upper}, and then, \cref{theorem:SDP-duality} guarantees the optimality of $(Z,P)$ and zero duality gap. 

From \cref{subsection:duality-analysis}, we know that
$P^\star$ is feasible to the dual SDP
\cref{eq:LQ-game_Upper_Kstar-gram-dual} and thus satisfies \cref{eq:KKT-upper-dual-feasibility}. This $P^\star$ is the dual candidate. We now construct a primal candidate  $Z^\star$. 
For \cref{eq:LQ-game_Upper_Kstar}, the disturbance $\nu (t; x(s),\,s\leq t)=w(t)= L^\star x(t)$ is  feasible due to the stability of $A+BK^\star + B_w L^\star$. 
Thus, 
the Gramian matrix
\begin{align*}
Z^\star = &\int_0^\infty
\begin{bmatrix}
    x(t)\\
    L^\star x(t)
\end{bmatrix}
\begin{bmatrix}
    x(t)\\
    L^\star x(t)
\end{bmatrix}^\tr dt
\begin{bmatrix}
    I\\
    L^\star
\end{bmatrix}Z_{11}^\star 
\begin{bmatrix}
    I\\
    L^\star
\end{bmatrix}^\tr
\end{align*}
with $Z_{11}^\star = \int_0^\infty x(t)x(t)^\tr dt$
is also feasible to the primal SDP \cref{eq:LQ-game_Upper_Kstar-gram} and thus satisfies \cref{eq:KKT-upper-primal-feasibility}.

We next verify that this primal and dual pair $Z^\star, P^\star$ satisfy the complementarity slackness \cref{eq:KKT-upper-complementary-slackness}. 
Since $K^\star = - R^{-1}B^\tr P^\star$ and $Q^\star = Q + \KstarTr R \Kstar$, we have
$$
\begin{aligned}
% &
\ABKstar^\tr P^\star+P^\star \ABKstar + Q^\star
=\underbrace{\mathcal{R}(P^\star)}_{=0} -\frac{1}{\gamma^2}P^\star B_wB_w^\tr P^\star.
\end{aligned}
$$ 
Thus, the left hand side of \cref{eq:KKT-upper-complementary-slackness} for $P^\star$ is reduced to
\begin{align}
&\left\langle Z^{\star},\begin{bmatrix}
        -\frac{1}{\gamma^2}P^\star B_wB_w^\tr P^\star & P^\star B_w\\
        B_w^\tr P^\star & -\gamma^2 I
    \end{bmatrix}\right\rangle \nonumber \\
    =& -\Tr\left(
    Z^\star
    \begin{bmatrix}
        \gamma \Lstar^\tr \\
        -\gamma I
    \end{bmatrix}
        \begin{bmatrix}
        \gamma \Lstar^\tr \\
        -\gamma I
    \end{bmatrix}^\tr
    \right), \nonumber \\
    =& -\Tr\Biggl(
    \begin{bmatrix}
    I\\
    L^\star
\end{bmatrix}Z_{11}^\star 
\underbrace{\begin{bmatrix}
    I\\
    L^\star
\end{bmatrix}^\tr
    \begin{bmatrix}
        \gamma \Lstar\\
        -\gamma I
    \end{bmatrix}}_{=0}
        \begin{bmatrix}
        \gamma \Lstar\\
        -\gamma I
    \end{bmatrix}^\tr
    \Biggr) =0,
    \label{eq:upper-KKT-cs-factorization}
\end{align}
where we have used the fact $\gamma \Lstar^\tr  = \frac{1}{\gamma}P^\star B_w$.

Hence,
from \cref{theorem:SDP-duality}, 
the pair $(Z^\star,P^\star)$ is optimal for the SDPs \cref{eq:LQ-game_Upper_Kstar-gram,eq:LQ-game_Upper_Kstar-gram-dual}, and strong duality holds with $\pol_\Kstar^\sdp = \overline{d}_\Kstar=x_0^\tr P^\star x_0$.
From the construction of $Z^\star$, we clearly have $\pol_\Kstar=\pol_\Kstar^\sdp$, and hence, by $\pol_\Kstar\leq \pol_\Kstar^\sdp$,
an optimal solution to \cref{eq:LQ-game_Upper_Kstar} is given by $\nu(t;x(\tau),\,\tau\leq t)=L^\star x(t)$. 
Consequently,
$\max_{\nu\in\CN}J_\gamma(\mu^\star,\nu) = x_0^\tr P^\star x_0=\pol^\star$ and a maximizer is 
$\nu^\star = L^\star x.$

The lower value $\pul^\star = \min_{\mu\in\CM} J_\gamma(\mu,\nu^\star)$ with \cref{eq:argmin=ustar} 
directly follows from the LQR analysis in \cref{subsection:gramian-LQR}. This is because \cref{eq:LQ-game_Lower_L-gram} becomes an LQR when $Q-\frac{1}{\gamma^2}P^\star B_wB_w^\tr \succ 0$.

\begin{remark}
In non-cooperative games with closed-loop information structures, 
the condition 
$Q-\frac{1}{\gamma^2}P^\star B_wB_w^\tr P^\star \succ 0$
is crucial for the saddle point property of the linear policy pair $(\mu^\star,\nu^\star)$.
In non-cooperative games, we should avoid an explicit condition  $\lim_{t\to\infty}x(t)\to0$, since it may require cooperation between the two players and the game is no longer fully non-cooperative. We also note that characterizing the upper value $\pol^\star$ does not require $Q-\frac{1}{\gamma^2}P^\star B_wB_w^\tr P^\star \succ 0$. Thus, we can extend our primal and dual analysis to $\Hinf$ control. This will be detailed in \cref{section:Hinf}. 
\hfill$\square$
\end{remark}

\section{Application to $\Hinf$ control}\label{section:Hinf}

In this section, we discuss 
$\Hinf$ control
as an extension of \cref{theorem:main_theorem} and the proof in \cref{section:proof-main-theorem}.
Here, we consider an zero initial condition $x(0)=0$ and define the output $z$ as the performance measure:
\begin{equation*}
    z(t) = \begin{bmatrix}
        Q^{1/2} &
        0
    \end{bmatrix}^\tr x(t) + \begin{bmatrix}
        0&
        R^{1/2}
    \end{bmatrix}^\tr u(t).
\end{equation*}

\subsection{$\Hinf$ control as an LQ game}
In control theory,
$\Hinf$ control is one of the most fundamental problems \cite{zhou1996robust,green2012linear}.  
This problem aims to synthesize a stabilizing controller that is robust against disturbances in the sense of
\begin{equation}\label{eq:Hinf_criterion}
    \|\bT_{zw}\|_{\Hinf} <\gamma
\end{equation}
for given $\gamma>0$. Here, $\bT_{zw}:\domwL\to\mathcal{L}_2^{n+m}[0,\infty)$ is the transfer function matrix from $w$ to $z$, and $\|\cdot\|_{\Hinf}$ denotes the \textit{$\Hinf$ norm}.   
It is well known that we have 
\begin{equation}\label{eq:Hinf-norm_def}
    \|\bT_{zw}\|_{\Hinf}
    =\sup_{
w\in\mathcal{L}_2^p[0,\infty)\setminus\{0\}}
\|z\|_2/\|w\|_2.
\end{equation}

The relationship \cref{eq:Hinf-norm_def} informs us that \cref{eq:Hinf_criterion} is equivalent to
\begin{equation*}
    \|z\|_2^2 -\gamma^2\|w\|_2^2 
=\int_0^\infty \left(
    x^\tr Qx+u^\tr Ru
    -\gamma^2\|w\|^2\right)dt   
    <0
\end{equation*}
for all $w\in\domwL,\, w\neq 0$.  
Thus, we can view $\Hinf$ control as finding a control $u$ that gives a negative $\pol^\star$ in \cref{eq:LQ-game-upper-value} for $w\neq 0$.
Accordingly, \cref{theorem:main_theorem} with $x_0=0$ immediately gives a sufficient condition for
\begin{equation*}
\|\bT_{zw}\|_{\Hinf}\leq \gamma.
\end{equation*}
However,
 establishing the strict inequality $\|\bT_{zw}\|_{\Hinf} <\gamma$ requires some additional analysis. 

We have the following theorem guaranteeing that
$u(t) = K^\star x(t)$ also serves as an $\Hinf$ controller solving \cref{eq:Hinf_criterion}.

\vspace{2pt}
\begin{theorem}\label{theorem:Hinf-control}
Consider the system \cref{eq:dynamics} with 
$x(0)=x_0=0$.
Assume that
the ARE in \cref{eq:indefinite-ARE-gamma} has a solution $P^\star\succeq0$ and $A-BR^{-1}B^\tr P^\star +\frac{1}{\gamma^2}B_wB_w^\tr P^\star$ is stable.
Then, 
for
$u(t)=K^\star x(t)= -R^{-1}B^\tr P^\star x(t)$,
we have
$\|\bT_{zw}\|_{\Hinf} <\gamma$.
\end{theorem}
\vspace{2pt}

The classical proof is based on the bounded real lemma \cite[Theorem 6.3.2]{green2012linear}. We here provide a new proof based on SDPs 
and the factorization \cref{eq:upper-KKT-cs-factorization}. %

\subsection{SDP-based proof for \cref{theorem:Hinf-control}}

Here, we fix the control input $u(t)=K^\star x(t)$, and thus we have 
$$z(t)
    = \begin{bmatrix}
        Q^{1/2},(R^{1/2}K^\star)^\tr
\end{bmatrix}^\tr x(t).$$
We follow the same Gramian approach in \Cref{subsection:SDP-upper-lower-values} and derive the SDP \cref{eq:LQ-game_Upper_Kstar-gram}. 
We only need to prove that the optimal value $\pol_{\Kstar}< 0$ when $x_0 = 0$ and $Z \neq 0$ in \cref{eq:LQ-game_Upper_Kstar-gram}. 
Note that $Z \neq 0$ corresponds to a nonzero disturbance $w$.

Recall that $P^\star$ is the stabilizing solution to the ARE \cref{eq:indefinite-ARE-gamma}. Then for any feasible solution $Z\succeq 0$ to \cref{eq:LQ-game_Upper_Kstar-gram}, we have 
\begin{align}
&\langle Z,\diag(Q^\star,-\gamma^2 I) \rangle \nonumber \\
=& \langle Z,\diag(Q^\star,-\gamma^2 I) \rangle \nonumber \\
&+\left\langle\underbrace{
\ABKstar Z_{11}+ B_wZ_{12}^\tr
+ \left(\ABKstar Z_{11}+ B_wZ_{12}^\tr\right)^\tr 
}_{=0},P^\star\right\rangle \nonumber \\
=& \left\langle 
Z,\begin{bmatrix}
\ABKstar^\tr  P^\star + P^\star \ABKstar 
        + Q^\star
        % +Q+\KstarTr R\Kstar
        & P^\star B_w\\
        B_w^\tr P^\star & - \gamma^2I
    \end{bmatrix}
    \right\rangle \nonumber \\
    =&- \left\langle 
Z,\begin{bmatrix}
        \gamma \Lstar^\tr \\
        -\gamma I
    \end{bmatrix}
        \begin{bmatrix}
        \gamma \Lstar^\tr \\
        -\gamma I
    \end{bmatrix}^\tr
    \right\rangle \leq 0. \label{eq:Hinf-proof-1}
\end{align}
This indicates that the optimal value $\pol_{\Kstar}^\sdp \leq 0$ in \cref{eq:LQ-game_Upper_Kstar-gram}. 

For any disturbance $w\in\domwL$, we set
\begin{equation} \label{eq:disturbance-deviation}
    \hat{w}(t)=w(t)-L^\star x(t), \,t\geq 0.
\end{equation} 
The stability of $\ABKstar + B_wL^\star $ ensures
$w\in\domwL\Leftrightarrow \hat w\in\domwL$.
The Gramian matrix $\hat{Z}$ below is feasible to \cref{eq:LQ-game_Upper_Kstar-gram}:
\begin{align*}
\hat{Z}
=& \int_0^\infty
\begin{bmatrix}
    x(t)\\
    w(t)
\end{bmatrix}
\begin{bmatrix}
    x(t)\\
    w(t)
\end{bmatrix}^\tr
dt \\
=& \int_0^\infty
\begin{bmatrix}
    x(t)\\
    L^\star x(t) +\hat{w}(t)
\end{bmatrix}
\begin{bmatrix}
    x(t)\\
    L^\star x(t) +\hat{w}(t)
\end{bmatrix}^\tr
dt \\
=& \begin{bmatrix}
    I\\
    L^\star
\end{bmatrix}Z_{11}\begin{bmatrix}
    I\\
    L^\star
\end{bmatrix}^\tr
+ \begin{bmatrix}
 0 & \hat{Z}_{12}   \\
 \hat{Z}_{12}^\tr &
 \Lstar\hat{Z}_{12} + \hat{Z}_{12}^\tr \Lstar^\tr
 % +\hat{Z}_{12}^\tr \LstarTr 
 +\hat{Z}_{22}   
\end{bmatrix},
\end{align*}
where $Z_{11} = \int_0^\infty x(t)x(t)^\tr dt$,
$\hat{Z}_{12} = \int_0^\infty x(t)\hat w(t)^\tr dt$,
and
$\hat{Z}_{22} = \int_0^\infty \hat w(t)\hat w(t)^\tr dt$.
Then, we observe that
\begin{align*}
&\langle \hat{Z},\diag(Q^\star,-\gamma^2 I) \rangle 
= \|z\|_2^2 - \gamma^2 \|w\|_2^2
\\
=&
- \left\langle 
\begin{bmatrix}
 0 & \hat{Z}_{12}   \\
 \hat{Z}_{12}^\tr &
 \Lstar\hat{Z}_{12} + \hat{Z}_{12}^\tr \Lstar^\tr
 % \Lstar\hat{Z}_{12}
 % +\hat{Z}_{12}^\tr \LstarTr 
 +\hat{Z}_{22}   
\end{bmatrix},\begin{bmatrix}
        \gamma \Lstar^\tr \\
        -\gamma I
    \end{bmatrix}
        \begin{bmatrix}
        \gamma \Lstar^\tr \\
        -\gamma I
    \end{bmatrix}^\tr
    \right\rangle\\
=& -\gamma^2 \Tr(\hat{Z}_{22})
= - \gamma^2\|\hat{w}\|^2_2,
\end{align*}
which indicates that
\begin{equation*}
\langle \hat{Z},\diag(Q^\star,-\gamma^2 I) \rangle 
=\|z\|_2^2 - \gamma^2 \|w\|_2^2
<0
\quad\Leftrightarrow
\quad\hat{w}\neq 0.
\end{equation*}

Since $x_0=0$, 
setting $\hat{w}=0$ implies $w = 0$; we thus have
$w\neq 0\Rightarrow  \hat{w}\neq 0$.
Therefore, we obtain
% \begin{equation*}
 $\|z\|_2^2 - \gamma^2 \|w\|_2^2<0$
% \end{equation*}
for any $w\in\domwL,\,w\neq 0$, which yields
$\|\bT_{zw}\|_{\Hinf}<\gamma$.

\begin{remark}[Relation with completion-of-squares]
A classical approach for computing an $\Hinf$ controller is the completion-of-squares technique.
The key trick is adding and subtracting
$\int_0^\infty\frac{d}{dt}x(t)^\tr P^\star x(t)$ 
for $\|z\|_2^2$.
In our approach, the computation of \cref{eq:Hinf-proof-1}, inspired by the Lagrange dual process, essentially plays the same role.
This implies that the completion-of-squares technique can be viewed as transforming a Lagrange function using $P^\star$.
\hfill$\square$
\end{remark}

\begin{remark}
In $\Hinf$ control, it is also important to show the converse statement of \cref{theorem:Hinf-control},
i.e., a stabilizing solution $P^\star$ to the ARE \cref{eq:indefinite-ARE-gamma} exists when \cref{eq:Hinf_criterion} is feasible.
It is possible to extend our approach to analyze this converse statement.
One potential way is to ensure the existence of a solution for primal and dual SDPs by utilizing the notion of \textit{strict feasibility}, which has a direct connection to the controllability, observability and initial conditions (see e.g., \cite[Section 3]{watanabe2025revisiting} and \cite[Proposition 5]{you2015primal}).
Due to the~page limit,
we leave the detailed investigation to our future work.~ 
\hfill$\square$
\end{remark}

\section{Illustrative example}\label{section:example}
\begin{figure*}[ht]
  \centering
  \begin{subfigure}{0.32\textwidth}
    \includegraphics[width=0.82\linewidth]{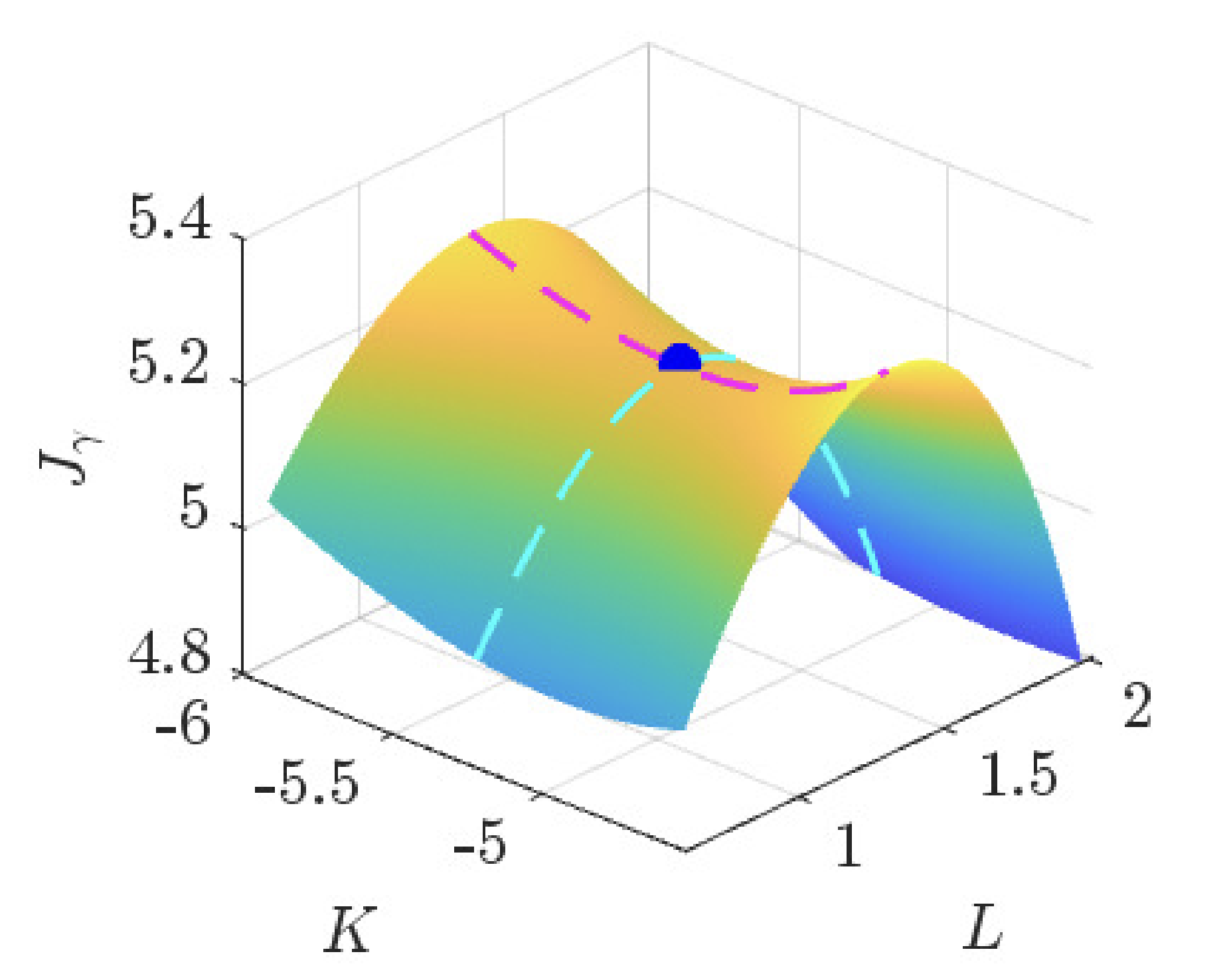}
    \caption{Saddle point: $q=10$ (enlarged)}
    \label{fig:sub1}
  \end{subfigure}
  \hfill
  \begin{subfigure}{0.32\textwidth}
    \includegraphics[width=0.85\linewidth]{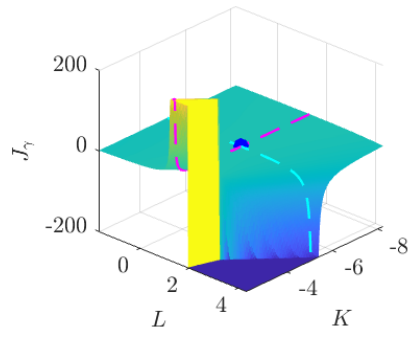}
    \caption{Saddle point: $q=10$}
    \label{fig:sub2}
  \end{subfigure}
  \hfill
  \begin{subfigure}{0.32\textwidth}
    \includegraphics[width=0.96\linewidth]{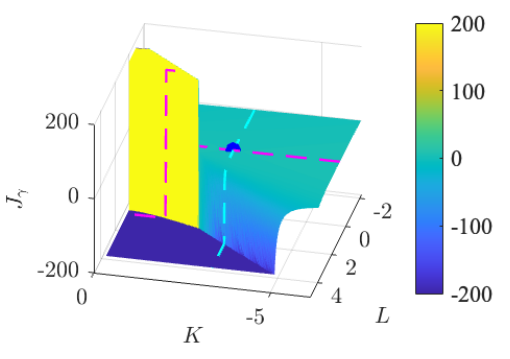}
    \caption{Non-saddle point: $q=1$}
    \label{fig:sub3}
  \end{subfigure}
  \caption{(a) The landscape of $J_\gamma$ 
  with $q=10$
  and the saddle point $(K^*,L^*)$ (denoted by the blue dot); (b) the landscape of $J_\gamma$ with $q=10$ over a broader region;
  (c) the landscape of $J_\gamma$ with $q=1$.
  The red and blue dashed lines represent $J_\gamma$ with $L=L^\star$ and $K=K^\star$, respectively. 
  Values are clipped to $200$ wherever $|J_\gamma|>200$.
  In (c), we can observe that $\min_{K}J_\gamma(Kx,L^\star x)=-\infty <J_\gamma(K^\star x,L^\star x)$ and thus $(K^*,L^*)$ does not constitute a saddle point.
  }
  \label{fig:landscape}
\end{figure*}

In deriving a saddle point of $J_\gamma$,
we have seen that
the assumption of $Q-\frac{1}{\gamma^2}P^\star B_wB_w^\tr P^\star\succ 0$ in \cref{theorem:main_theorem} is important. 
To illustrate the role of this assumption,
this section provides visualizations of the landscape of $J_\gamma$ with different parameters.

Consider arguably the simplest scalar case 
\[
\dot{x}=x+u+w,\qquad x(0)=1
\]
    and the objective function
    \[
        J_\gamma = \int_0^\infty (qx^2+u^2-\gamma^2w^2)dt,\quad \gamma=2.
    \]
We now consider two cases: $q=10$ and $1$, where $Q-\frac{1}{\gamma^2}P^\star B_wB_w^\tr P^\star\succ 0$ is satisfied and unsatisfied respectively.
Now, 
for $q>0$, the ARE in \cref{eq:indefinite-ARE-gamma} and the stabilizing solution $p^*$ are given by
        $3p^2-8p-4q=0$ and $p^\star = \frac{4}{3}+\frac{2}{3}\sqrt{4+3q}$.
    We thus obtain $(\mu^\star, \nu^\star)=\left(K^\star x,L^\star x\right)$ with
    \[
    (K^\star,L^\star) = \left(-p^*, p^*/4\right).
    \]
    Note that when restricting $u$ and $w$ to be linear as $u=Kx$ and $w=Lx$ with $K,L\in\mathbb{R}$, we can explicitly write 
    \begin{align*}
    J_\gamma =
    \begin{cases}
       (q + K^2-\gamma^2 L^2) \times \infty, &
       1+K+L\geq 0\\
       \displaystyle\frac{q + K^2-\gamma^2 L^2 }{2|1+K+L|}, 
       &1+K+L< 0.
    \end{cases}
    \end{align*}

    We plot
    the landscapes of $J_\gamma$ for $(K,L)$
    in \cref{fig:landscape}.
    Here, \cref{fig:sub1,fig:sub2} show the case of $q=10$ and \cref{fig:sub3} corresponds to $q=1$. The blue dots represent $(K^\star,L^\star)$, and values are clipped to 200 whenever $|J_\gamma|>200$.
    The red and blue dashed lines represent the values of $J_\gamma$ with $L=L^\star$ and $K=K^\star$, respectively. 
    % In \cref{fig:sub1} and \cref{fig:sub3}, 
    % the landscape is $q=200$ is plotted, and \cref{fig:sub3} 
    From \cref{fig:sub1,fig:sub2}, we observe that $(K^*,L^*)$ is a saddle point for $q=10$, as expected from \cref{theorem:main_theorem}. In contrast, for $q=1$, $(K^\star,L^\star)$ does not constitute a saddle point (\cref{fig:sub3})
    because for $L=L^\star$, the control player can enforce
    $J_\gamma=-\infty$ by choosing $K = 0$  (see the red line).
    We note that
    even in the case of $q=1$, $L=L^\star$ still maximizes $J_\gamma(K^\star,L)$ as the blue~line.
    These results show the importance of $Q-\frac{1}{\gamma^2}P^\star B_wB_w^\tr P^\star\succ 0$ for the saddle point property of the linear policy pair \cref{eq:u-w-pair}.

\section{Conclusion}\label{section:conclusion}
This paper addressed the infinite-horizon LQ differential game for continuous-time systems 
from the perspective of SDP duality.
Under the assumption that a stabilizing solution to a Riccati equation exists with a regularity condition,
we presented a novel SDP duality-based proof for
the saddle point property given by linear static policies.
In the proof,
we leveraged the Gramian representation technique and clarified
the role of the ARE, dual problems, and KKT condition.
Moreover, we applied these results to
the $\Hinf$ suboptimal control problem.
Our future directions include
a more detailed analysis of the LQ game and $\Hinf$ control,
such as 
the existence of the stabilizing solution to the ARE.
\bibliographystyle{IEEEtran}
\bibliography{ref.bib}

\end{document}